\documentclass[12pt]{smfart}
\usepackage{color}
\usepackage{amssymb,verbatim}
\usepackage{hyperref}
\usepackage{amsthm,array,amssymb,amscd,amsfonts,amssymb,latexsym, url}
\usepackage{amsmath}
\usepackage[all]{xy}
\usepackage[french]{babel}

\newcounter{spec}
{\end{list}}

\newcommand{\A}{{\mathbf A}}
\renewcommand{\P}{{\mathbf P}}

\renewcommand{\L}{{\mathbb L}}
\newcommand{\oi}{\hskip1mm {\buildrel \simeq \over \rightarrow} \hskip1mm}
\newcommand{\Br}{{\operatorname{Br} }}

\renewcommand{\lim}{\varprojlim}

\numberwithin{equation}{section}

\newfont{\gothic}{eufb10}

\renewcommand{\qed}{{\hfill$\square$}}

\newtheorem{theo}{Th\'{e}or\`{e}me}[section]
\newtheorem{prop}[theo]{Proposition}

\newtheorem{lem}[theo]{Lemme}
\newtheorem{cor}[theo]{Corollaire}
\theoremstyle{definition}
\newtheorem{defi}[theo]{Definition}
\theoremstyle{remark}

\newcommand{\bthe}{\begin{theo}}
\newcommand{\ble}{\begin{lem}}
\newcommand{\bpr}{\begin{prop}}
\newcommand{\bco}{\begin{cor}}
\newcommand{\bde}{\begin{defi}}
\newcommand{\ethe}{\end{theo}}
\newcommand{\ele}{\end{lem}}
\newcommand{\epr}{\end{prop}}
\newcommand{\eco}{\end{cor}}
\newcommand{\ede}{\end{defi}}

\DeclareFontFamily{U}{wncy}{}
\DeclareFontShape{U}{wncy}{m}{n}{%
<5>wncyr5%
<6>wncyr6%
<7>wncyr7%
<8>wncyr8%
<9>wncyr9%
<10>wncyr10%
<11>wncyr10%
<12>wncyr6%
<14>wncyr7%
<17>wncyr8%
<20>wncyr10%
<25>wncyr10}{}
\DeclareMathAlphabet{\cyr}{U}{wncy}{m}{n}

\begin{document}
 
 \title[Rationalit\'e d'un fibr\'e en coniques]{Rationalit\'e d'un fibr\'e en  coniques}
 
\author{Jean-Louis Colliot-Th\'el\`ene}
\address{CNRS, Universit\'e Paris-Sud\\Math\'ematiques, B\^atiment 425\\91405 Orsay Cedex\\France}
\email{jlct@math.u-psud.fr}

\date{arXiv v1, 21 octobre 2013; pr\'esente version, 21 novembre 2013}

\maketitle

\section{Introduction}

Soit $k$ un corps de caract\'eristique diff\'erente de 2
 contenant un \'el\'ement $i$ de carr\'e $-1$.
 La cubique $C$  d'\'equation affine 
$$y^2=x(x^2-1)$$
sur le corps $k$ admet l'automorphisme $g$, d'ordre 4,  d\'efini par  $$(x,y) \mapsto (-x,iy).$$
On consid\`ere le quotient du produit triple $C\times_{k}C\times_{k}C$ 
par l'action diagonale de $g$.

Motiv\'e par un travail ant\'erieur de Ueno \cite{U}, Campana \cite{C} montre, pour $k$ le corps des complexes, que 
ce quotient est 
 une vari\'et\'e rationnellement connexe et demande si cette vari\'et\'e
est rationnelle. Pour $k$ arbitraire,
Catanese, Oguiso et  Truong \cite{COT} (arXiv, 14 oct. 2013) montrent ensuite que  c'est une vari\'et\'e
$k$-unirationnelle. Pour ce faire, ils commencent par montrer \cite[Prop. 2.5]{COT} que
le quotient est $k$-birationnel   \`a  la $k$-vari\'et\'e $X$ d\'efinie dans l'espace affine
$\A^4_{k}$ avec  coordonn\'ees $s,t,u,v$ par l'\'equation
$$s^2 v(v^2-1) -t^2 u(u^2-1) + uv(u^2-v^2)=0.$$

Dans cette note,  je montre que la $k$-vari\'et\'e $X$  est $k$-rationnelle,
c'est-\`a-dire que son corps des fonctions est transcendant pur sur $k$.
Ceci r\'epond \`a la question de Campana.

La rationalit\'e de ce genre de quotient a r\'ecemment attir\'e l'attention
des sp\'ecialistes de dynamique complexe   \cite{OT}.

La $k$-vari\'et\'e $X$ est fibr\'ee en coniques sur le  plan affine de coordonn\'ees $u,v$.
Fr\'ed\'eric Campana avait sugg\'er\'e de calculer la cohomologie non
ramifi\'ee de $X$, 
et en particulier son groupe de Brauer non ramifi\'e (pour un rapport de synth\`ese sur la cohomologie non ramifi\'ee, voir \cite{CT}). Pour les fibrations en coniques,
  on peut suivre
  la m\'ethode d'Ojanguren et l'auteur \cite{CTO}.
Lors de ce calcul, 
 il est  apparu que la fibration en coniques concern\'ee
est $k$-birationnellement \'equivalente, au-dessus du plan affine, \`a une autre
fibration en coniques dont l'espace total est clairement $k$-rationnel.
  
 Je remercie Fr\'ed\'eric Campana de m'avoir inform\'e de ce probl\`eme
 et de ses d\'eveloppements.

\section{Rappels}

Donnons  quelques rappels sur les coniques et sur le groupe de Brauer des corps
et des sch\'emas. On renvoie le   lecteur 
aux ouvrages  \cite{GB}, \cite{S}, \cite{GSz}.

On note $G[2]$ le sous-groupe des \'el\'ements de 2-torsion
d'un groupe ab\'elien $G$.

\medskip

Soit $F$ un corps de caract\'eristique diff\'erente de 2. 
\`A  la donn\'ee de deux \'el\'ements $a,b \in F^{\times}$ on associe
la classe de l'alg\`ebre de quaternions $(a,b) $ dans 
le groupe de Brauer $ \Br F$. C'est un \'el\'ement de 2-torsion.

On note $C_{a,b} \subset \P^2_{F}$ la conique
d'\'equation projective 
$$ X^2-aY^2-bT^2=0.$$

La classe $(a,b)$ est nulle dans  $ \Br(F)$ si et seulement si la conique $C_{a,b}$
sur le corps $F$  
a un point rationnel, ce qui est encore \'equivalent au fait que
son corps des fonctions est transcendant pur sur $F$.

Plus g\'en\'eralement, on a la proposition suivante.

\begin{prop}\label{coniques}
Avec les notations ci-dessus, \'etant donn\'es quatre  \'el\'ements $a,b, a',b' \in F^{\times}$,
les \'enonc\'es suivants sont \'equivalents :

(i) $(a,b)=(a',b') \in \Br(F)$.

(ii) Les coniques $C_{a,b}$ et $C_{a',b'}$ sont isomorphes sur le corps $F$.

(iii) Il existe une homographie, c'est-\`a-dire un \'el\'ement de ${\rm PGL}_{3}(F)$
qui transforme la conique $C_{a,b} \subset \P^2_{F}$ en la conique $C_{a',b'} \subset \P^2_{F}$.
\end{prop}

L'\'equivalence de (i) et (ii) remonte \`a Witt \cite[\S 2]{W}, voir \cite[Thm. 1.4.2]{GSz}.
Soit $f :  C_{1 } \oi C_{2}$ un $k$-isomorphisme abstrait de coniques lisses sur le corps $k$.
L' inverse du faisceau canonique $\omega_{2}$ sur $C_{2}$ est le faisceau canonique $\omega_{1}$ sur $C_{1}$.
Tout plongement  de $C_i$ dans  $\P^2_{k}$ comme conique est associ\'e au choix d'une base de l'espace vectoriel $H^0(C_{i},\omega_{i}^{-1})$, qui est de dimension 3.
Ceci montre que (ii) implique (iii).

\medskip

Soit $A$ un anneau de valuation discr\`ete de corps des fractions $K$,
de corps r\'esiduel $\kappa$ de caract\'eristique diff\'erente de 2.
On dispose d'une application r\'esidu
$$\partial_{A} :  \Br (K) [2] \to \kappa^{\times}/\kappa^{\times 2}$$
qui envoie une classe de quaternions $(a,b)$ (avec $a, b \in K^{\times}$)
sur la classe de $$(-1)^{v(a)v(b)}  [a^{v (b)}/b^{v(a)}],$$
o\`u $[a^{v (b)}/b^{v(a)}]$ d\'esigne la classe dans $\kappa^{\times}/\kappa^{\times 2}$ de l'unit\'e
$a^{v (b)}/b^{v(a)} \in A^{\times}$. L'application r\'esidu s'ins\`ere dans
une suite exacte
$$0 \to \Br(A)[2]  \to \Br(K)[2]  \to \kappa^{\times}/\kappa^{\times 2} \to 1.$$

\medskip

Soit $X$ une $k$-vari\'et\'e lisse int\`egre de corps des fonctions $k(X)$.
\`A chaque point $x$ de codimension 1 on associe son anneau local,
qui est un anneau de valuation discr\`ete de corps des fractions $k(X)$,
de corps r\'esiduel $\kappa(x)$, le corps des fonctions de la sous-vari\'et\'e
d\'efinie par $x$.

Supposant ${\rm car}(k) \neq 2$, on a une suite exacte
$$0 \to \Br (X)[2] \to \Br ( k(X))[2]  \to  \bigoplus_{x \in X^{(1)}} \kappa(x)^{\times}/\kappa(x)^{\times 2},$$
o\`u la fl\`eche de droite est celle d\'efinie par les divers r\'esidus en les points $x \in X^{(1)}$.

\medskip

Pour $X= \A^n_{k}$ l'espace affine sur un corps de caract\'eristique diff\'erente de 2,
  l'application naturelle $\Br(k)[2] \to \Br(\A^n_{k})[2]$ est un isomorphisme.

\medskip 

On a donc :

\begin{prop}\label{espaceaffine}
Soit $k$ un corps de caract\'eristique diff\'erente de 2, et soit $n\geq 1$ un entier.
Les applications r\'esidus aux points de codimension 1 de l'espace affine $\A^n_{k} $ donnent naissance
\`a la suite exacte :
$$0 \to \Br(k)[2] \to \Br(k(\A^n))[2] \to  \bigoplus_{x \in X^{(1)}} \kappa(x)^{\times}/\kappa(x)^{\times 2},$$
\end{prop}

On peut donner de ce r\'esultat une formulation et une d\'emonstration plus \'el\'ementaire
en utilisant   la cohomologie galoisienne \cite{S}. Le cas $n=1$ est un cas
particulier de la suite exacte de Faddeev \cite[Cor. 6.4.6]{GSz}. Le cas $n>1$
peut se ramener au cas $n=1$ 
par fibrations successives sur des espaces affines de dimension plus petite.

\section{Le calcul}

On consid\`ere la conique sur  le corps $k(u,v)$ d\'efinie par
l'\'equation affine (variables $s,t$)
$$s^2 v(v^2-1) -t^2 u(u^2-1) + uv(u^2-v^2)=0$$
soit encore
$$s^2-uv(u^2-1)(v^2-1) t^2 - u(v^2-1)(v^2-u^2)=0.$$

On associe \`a cette conique la classe de quaternions
$$\alpha = (uv(u^2-1)(v^2-1), u (v^2-1)(v^2-u^2)) \in \Br  (k(u,v)).$$

Les seuls r\'esidus non triviaux possibles aux points de codi\-mension 1
de $\A^2_{k}$ sont aux points donn\'es par l'annulation de l'un des
facteurs $uv(u^2-1)(v^2-1)$  ou $u (v^2-1)(v^2-u^2)$.

Appliquant la formule rappel\'ee au paragraphe 2, on trouve  les r\'esidus suivants

En $u=0$, la classe de $v $ dans $k(v)^{\times}/k(v)^{\times 2}$.

En $v=0$, la classe de $u$ dans  $k(u)^{\times}/k(u)^{\times 2}$.

En $v-1=0$, la classe $1$ dans $k(u)^{\times}/k(u)^{\times 2}$.

En $v+1=0$, la classe $-1$ dans $k(u)^{\times}/k(u)^{\times 2}$.

En $u-1=0$, la classe $1$ dans $k(v)^{\times}/k(v)^{\times 2}$.

En $u+1=0$, la classe $-1$ dans $k(v)^{\times}/k(v)^{\times 2}$.

En $u-v=0$, la classe $1$ dans $k(u)^{\times}/k(u)^{\times 2}$.

En $u+v=0$, la classe $-1$ dans $k(u)^{\times}/k(u)^{\times 2}$.

Par ailleurs,  les seuls  r\'esidus de    $(u,v) \in \Br(k(u,v))$ 
aux points de codimension 1 de $\A^2_{k}$ qui ne sont pas \'egaux \`a $1$
sont 
 la classe de $v$ dans $k(v)^{\times}/k(v)^{\times 2}$ en $u=0$ et la classe de $u$ 
 dans 
 $k(u)^{\times}/k(u)^{\times 2}$ en $v=0$.
 
 Comme on a suppos\'e que  $-1=1 \in k^{\times}/k^{\times 2}$,
 on voit que l'\'el\'ement $$\alpha - (u,v) \in \Br ( k(u,v))$$
a tous ses r\'esidus
triviaux aux points de codimension 1 de $\A^2_{k} $, donc par la  Proposition \ref{espaceaffine}
 est la classe d'un \'el\'ement de $\gamma \in \Br (k)$. 
 
Pour calculer $\gamma$, on consid\`ere la restriction de  $(u,v)$ et de  
  $$(uv(u^2-1)(v^2-1), u (v^2-1)(v^2-u^2))$$
  au point g\'en\'erique de la droite $D$ d'\'equation $v=1-u$,
de corps des fonctions $k(u)$.  On a 
$$(u,1-u)=0 \in \Br (k(u)).$$
La restriction de
 $(uv(u^2-1)(v^2-1), u (v^2-1)(v^2-u^2)) $ au point g\'en\'erique de la droite $D$
  est \'egale \`a
 $$(u(1-u)(u^2-1)(u^2-2u), u(u^2-2u)(1-2u) ) = (-(u+1)(u-2),(u-2)(1-2u)) $$
 dans  $ \Br( k(u))$.

La restriction de  $\gamma= \alpha-(u,v)$ \`a la droite  $D$ 
a donc pour image   $(-(u+1)(u-2),(u-2)(1-2u)) \in \Br (k(u))$.
Ce dernier symbole est dans le groupe de Brauer de l'anneau de valuation discr\`ete
anneau local 
de la droite $D$  au point $P$ de coordonn\'ees $u=0, v=1$. Ainsi $\gamma$ et  $(-(u+1)(u-2),(u-2)(1-2u))$
sont \'egaux dans $\Br(A)$. 
 En \'evaluant au point $P$,
on trouve
 $$\gamma = (2, -2) =0.$$

On a donc $\alpha=(u,v) \in \Br (k(u,v))$.

\medskip

{\em Remarque.} 
 Sur un corps $k$ avec $\Br(k)=0$, tel que le corps des complexes ou plus g\'en\'eralement
un corps de dimension cohomologique 1, le calcul ci-dessus n'est pas n\'ecessaire.
Sur un corps quelconque de caract\'eristique diff\'erente de 2,3,5, avec $i \in k$, on peut directement
\'evaluer $\alpha$ au point $k$-rationnel $M$ de coordonn\'ees $u=2, v=3$.
Cela donne
$\alpha(M)=(2.3.3.8, 2.5.8)=0 \in \Br(k)$ puisque $2.3.3.8$ est un carr\'e.
Par ailleurs $(u,v)(M)=(2,3)=0$ car la conique d'\'equation $$X^2-2Y^2-3T^2=0$$
poss\`ede le point rationnel  $(X,Y,T)=(1,i,1)$.
\medskip

D'apr\`es la proposition \ref{coniques},  l'\'egalit\'e  $\alpha=(u,v) \in \Br (k(u,v))$
implique que
 les coniques dans $\P^2_{k(u,v)}$ donn\'ees en coordonn\'ees homog\`enes $S,T,R$ par
$$ v(v^2-1) S^2 - u(u^2-1) T^2 + uv(u^2-v^2) R^2=0$$
et
$$S^2-uT^2-vR^2=0$$
sont   isomorphes sur le corps $k(u,v)$.

En termes g\'eom\'etriques, cela signifie qu'au-dessus du plan affine $\A^2_{k}$
 de coordonn\'ees $u,v$, 
 les deux fibrations en coniques
 associ\'ees
sont birationnelles par un isomorphisme qui pr\'eserve la fibration.

La $k$-vari\'et\'e d'\'equation affine
$$s^2-ut^2-v=1$$
est clairement $k$-isomorphe \`a l'espace affine $\A^3_{k}$
de coordonn\'ees $s,u,t$.  On en conclut que la $k$-vari\'et\'e
d'\'equation affine
 $$s^2 v(v^2-1) -t^2 u(u^2-1) + uv(u^2-v^2)=0$$
 est  $k$-rationnelle.  D'apr\`es \cite[Prop. 2.5]{COT},
 ceci implique que
 le quotient de $C\times_{k}C\times_{k}C$ par l'automorphisme diagonal $g$
 d'ordre 4
 consid\'er\'e au-d\'ebut de cette note est une $k$-vari\'et\'e $k$-rationnelle.

\noindent
CNRS, UMR 8628, Math\'ematiques, B\^atiment 425, Universit\'e Paris-Sud,
F-91405 Orsay, France
\medskip

\noindent jlct@math.u-psud.fr

\end{document}